\documentclass[leqno, 12pt]{amsart}
\usepackage{amsmath, amscd, amssymb, mathrsfs, float}

\numberwithin{equation}{section}
\theoremstyle{plain}
  \newtheorem{thm}{Theorem}[section]
  \newtheorem{lem}[thm]{Lemma}
  \newtheorem{cor}[thm]{Corollary}
  \newtheorem{prop}[thm]{Proposition}

\theoremstyle{definition}
  \newtheorem{defn}[thm]{Definition}

\theoremstyle{remark}
  \newtheorem{rem}[thm]{Remark}
  \newtheorem{pf}{Proof}

\newcommand{\A}{\mathcal{A}}
\renewcommand{\b}{{\frak b}}
\newcommand{\C}{{\mathbb C}}
\newcommand{\CH}{{CH^{\cdot}}}
\renewcommand{\d}{{\frak d}}
\newcommand{\G}{{\mathbb G}}

\newcommand{\Z}{{\mathbb Z}}

\begin{document}

\title[]{The Chow ring of the moduli space\\
and its related homogeneous space\\
of bundles on ${\mathbb P}^2$ with charge $\text{\bf 1}$}
\author []{Yasuhiko Kamiyama}
\address[]{Department of Mathematics, University of the Ryukyus, 
Nishihara-Cho, Okinawa 903-0213, Japan} 

\author []{Michishige Tezuka}
\address[]{Department of Mathematics, University of the Ryukyus, 
Nishihara-Cho, Okinawa 903-0213, Japan} 
\email{tez@sci.u-ryukyu.ac.jp}

\subjclass[2000]{14M17 (14N10)} 
\keywords{Moduli space, homogeneous space, Chow ring, cycle map}
\maketitle
\begin{abstract} For an algebraically closed field $K$ with ${\rm ch} (K) \not = 2$, 
let ${\mathcal O} M(1,SO(n,K))$ denote the moduli space of holomorphic 
bundles on ${\mathbb P}^2$ 
with the structure group $SO(n,K)$ and half the first Pontryagin index being 
equal to $1$, 
each of which is trivial on a fixed line 
$l_\infty$ and has a fixed holomorphic trivialization there. In this paper we 
determine the Chow ring of ${\mathcal O} M(1,SO(n,K))$. 
\end{abstract}

\section{Introduction} Let $G$ be one of the classical groups 
$SU(n)$, $SO(n)$ or $Sp(n)$, and let $k \geq 0$ be half 
the first Pontryagin index of a $G$-bundle $P$ over 
$S^4= {\mathbb R}^4 \cup \{ \infty \}$. 
Denote by $M(k, G)$ the framed moduli space whose points represent 
isomorphism classes of pairs: 
\begin{equation*}
\text{(self-dual $G$-connections on $P$, isomorphism $P_\infty \simeq G$).}
\end{equation*}
Let ${\mathcal O} M(k, G^{\C})$ denote the moduli space of holomorphic 
bundles on ${\mathbb{CP}}^2$ for the associated complex group, trivial on a fixed line 
$l_\infty$ and with a fixed holomorphic trivialization there. Then Donaldson 
(\cite{Do}) showed a diffeomorphism $M (k, G) \simeq {\mathcal O} M(k,G^{\C})$. 
\par
In \cite{KKT} the topology of $M(1, SO(n)) \simeq {\mathcal O} M(1,SO(n, \C))$
was studied in detail. The result was used in \cite{K} to prove the fact that 
the natural homomorphism $J: H_\cdot (M(1, SO(n)), \Z/2) \to 
H_\cdot (\Omega_0^3 {\rm Spin} (n), \Z/2)$ is injective. 
Moreover, the image of $J$ was determined. 
To prove this, the following description of ${\mathcal O} M(1,SO(n, \C))$
by a homogeneous space was used: We set 
\begin{equation*}
W_n = SO(n)/ ( SO(n-4) \times SU(2)). 
\end{equation*}
Then there is a diffeomorphism 
\begin{equation}
{\mathcal O} M(1,SO(n, \C)) \simeq {\mathbb R}^5 \times W_n.
\end{equation}
\par
The purpose of this paper is to generalize the definition of 
${\mathcal O} M(1,SO(n, \C))$ for any algebraically closed field $K$ with 
${\rm ch} (K) \not = 2$ and to determine the Chow ring of this. 
The Chow ring of a classifying space was studied by Totaro \cite{To}. 
A loop space is considered to be a dual situation of a classifying space 
in a certain sense. Our result and the result of \cite{K} are the first step 
for a loop space. 
\begin{defn} Let $K$ be an algebraically closed field with ${\rm ch} (K) \not = 2$. 
Let ${\mathcal O} M(1,SO(n,K))$ denote the moduli space of holomorphic 
bundles on ${\mathbb P}^2$ 
with the structure group $SO(n,K)$ and half the first Pontryagin index being 
equal to $1$, 
each of which is trivial on a fixed line 
$l_\infty$ and has a fixed holomorphic trivialization there. 
\end{defn}
The moduli space ${\mathcal O} M(1,SO(n,K))$ is a quasi-projective variety and defines the Chow 
ring. More explicitly, the diffeomorphism (1.1) is generalized (in the sense of 
a biregular map) as follows: 
We set 
\begin{equation*}
X_n = SO (n, K)/ (SO (n-4,K) \times SL (2, K) ) \cdot P_u, 
\end{equation*}
where $P_u$ denotes the unipotent radical. (Recall that for a parabolic subgroup $P$ 
of an algebraic group $G$, $P$ is a semidirect product of a reductive group and
its unipotent radical $P_u$.) 
Then there is a biregular map 
\begin{equation}
{\mathcal O} M(1,SO(n, K)) \simeq {\mathbb A}^2 \times X_n. 
\end{equation}
(For the proof of (1.2), see Proposition 2.1.) A formula of Grothendieck \cite{C} 
shows that 
\begin{equation*}
\CH ({\mathcal O} M(1, SO(n,K))) \simeq \CH (X_n).
\end{equation*}
The purpose of this paper is to determine the Chow rings of $X_n$ 
and its related algebraic variety $Y_n$ explicitly.
\par
The Schubert cell approach of the Chow ring of $Y_n$ by using a Young diagram is done in \cite{PR1}, \cite{PR2}.
However it needs further work to determine the Chow ring of $X_n$ from this. Hence we first calculate the Chow ring 
of $Y_n$ more explicitly by a different method. Then we calculate the Chow ring from the results. Our results for the 
Chow ring of $X_n$ are new.
\par
This paper is organized as follows. In Sect.~2 we first prove (1.2). Then 
we recall basic facts on the Chow ring. 
In Sect.~3 we determine an integral basis and the ring structure of 
$\CH (Y_n)$, where $Y_n$ is an algebraic variety which is related to $X_n$. 
(See Theorems 3.7 and 3.9.)
The ring structure of $\CH (Y_n)$ proved in Theorem 3.9 is one of our main results.
Since the results are long, we give them in tables in Sect.~5. (See 5.2-5.5.)
Using the results of Sect.~3, we determine $\CH (X_n)$ in Sect.~4. (See Theorem 4.1.)
\par
We thank N.~Yagita for turning our interest to the Chow ring and explaining 
the paper \cite{SY}. 

\section{Preliminaries}
We first prove (1.2): 
\begin{prop} For an algebraically closed field $K$ with ${\rm ch} (K) \not = 2$, 
there is a biregular map 
\begin{equation*}
{\mathcal O} M(1, SO(n, K)) \simeq {\mathbb A}^2 \times X_n.
\end{equation*}
\end{prop}
\begin{pf} 
Recall that a monad description of ${\mathcal O} M(k, SO(n, \C))$ was
indicated in \cite{Do} and given explicitly in \cite{NS} and \cite{Ti}. It is easy 
to see that the description remains valid for any algebraically closed field $K$. 
In particular, ${\mathcal O} M(1, SO(n, K))$ is given as follows: 
\begin{lem}
Let $\mathscr{C}_n$ be the space of $n \times 2$ matrices 
\begin{equation*}
c= \begin{pmatrix}
z_1 & w_1\\
z_2 & w_2\\
\vdots & \vdots\\
z_n & w_n
\end{pmatrix}
\end{equation*}
with coefficients in $K$ satisfying: 
\begin{enumerate}
\renewcommand{\labelenumi}{\alph{enumi})}
\item $c^T c=O$, that is: 
\begin{equation*}
\sum_{i=1}^n z_i^2 =0, \quad
\sum_{i=1}^n w_i^2 =0
\quad \text{and} \quad 
\sum_{i=1}^n z_i w_i =0, 
\end{equation*}
\item The rank of $c$ over $K$ is $2$. 
\end{enumerate}
The group $SL(2, K)$ acts on $\mathscr{C}_n$ from the right by the multiplication of matrices. 
Then there is a biregular map
\begin{equation*}
{\mathcal O} M(1, SO(n, K)) \simeq {\mathbb A}^2 \times (\mathscr{C}_n/SL (2, K)). 
\end{equation*}
\end{lem}
From the lemma, it suffices to prove $X_n \simeq \mathscr{C}_n/SL (2, K)$. 
We prove this for the case $n= 2m$. (The case $n= 2m+1$ can be proved similarly.) 
Recall that in \cite{B}, $SO(n,K)$ was defined as follows: 
Let $q (x)$ be a quadratic form on ${\mathbb A}^n$ defined by 
$q(x) = \sum_{i=1}^m x_i x_{m+i}$, and let $B(x,y)$ be the associated bilinear form. 
Then $SO(n,K)$ is defined by 
\begin{equation*}
SO(n,K) = \{ \sigma \in \text{\rm Aut} ({\mathbb A}^n): 
B (\sigma(x), \sigma(y)) = B(x,y) \; \text{for} \; x, y \in {\mathbb A}^n \}.
\end{equation*}
We set 
\begin{equation*}
x_j= z_j+ \sqrt{-1} z_j, \; x_{m+j}= z_j- \sqrt{-1} z_j, \; y_j= w_j+ \sqrt{-1} w_j 
\quad \text{and} \quad
y_{m+j}= w_j-\sqrt{-1} w_j, 
\end{equation*}
where $1 \leq j \leq m$.
Then the defining equations of $\mathscr{C}_n$ are given by 
\begin{equation*}
q (x)= q(y)=0 \quad \text{and} \quad B(x,y) =0.
\end{equation*}
Clearly $SO(n,K)$ acts on $\mathscr{C}_n$. It is easy to prove the following lemma. 
(See \cite[V $23.4$]{B}.)
\begin{lem}
\begin{equation*}
SO(n,K)/ SO(n-4,K) \cdot P_u \simeq \mathscr{C}_n,
\end{equation*}
where $P_u$ is the unipotent radical of a parabolic subgroup with a Levi factor 
$SO(n-4, K) \times GL(2, K)$.
\end{lem}
Now Proposition 2.1 follows from Lemma 2.3. This completes the proof of Proposition 2.1. \qed
\end{pf}

Next we recall basic facts on the Chow ring. We suppose that 
an algebraic variety $V$ is defined over $K$. 
Let $\CH (V)$ denote the Chow ring and $CH^i (V)$ the subgroup of 
$\CH (V)$ generated by the cycles of codimension $i$. 

\begin{thm}[\cite{C}] {\rm (i)} Let $V$ be a nonsingular variety,  
$X$ a nonsingular closed subvariety of $V$, and $U = X-V$. Then 
there exists an exact sequence 
\begin{equation*}
\CH (X) \overset{i_\ast}{\to} \CH (V) \overset{j^\ast}{\to} 
\CH (U) \to 0,
\end{equation*}
where $i: X \to V$ (resp. $j: U \to V$) is a closed immersion 
(resp. an open immersion).
\par
For the definitions of $i_\ast$ and $j^\ast$, see also \cite{H}. 
\par
{\rm (ii)} Let $\pi: E \to V$ be a fiber bundle with an affine space 
${\mathbb A}^n$ as a fiber. Then the induced map $\pi^\ast: \CH (V) \to \CH (E)$ 
is an isomorphism. 
\end{thm}

The Chow ring of the following projective variety is well-known. 
\begin{thm}[\cite{BGG}, \cite{Dem}] Let $G$ be a reductive algebraic 
group and $P$ a maximal parabolic subgroup. Then 
\par
{\rm (i)} a quotient $G/P$ is a nonsingular projective variety. 
\par
{\rm (ii)} $\CH (G/P)$ is generated by the Schubert varieties. 
\par
{\rm (iii)} $\CH (G/P)$ is independent of ${\rm ch} (K)$. Moreover, 
$\CH (G/P) \simeq H^\cdot (G/P, {\mathbb Z})$ for $K= {\mathbb C}$. 

\end{thm}

\section{The ring structure of $\CH (Y_n)$} 
Before describing the results, we need some notations and results.
We set 
\begin{equation*}
Y_n = SO (n, K)/ (SO (n-4,K) \times GL (2, K) ) \cdot P_u. 
\end{equation*}
Then we have a principal bundle 
\begin{equation}
\G_m \to X_n \overset{\pi}{\to} Y_n.
\end{equation}
In this section we determine an integral basis and the ring structure of $\CH (Y_n)$.
By Theorem 2.5 (ii), (iii), we obtain the following theorem: 
\par

\begin{thm}[\cite{KKT}] We have an isomorphism as modules{\rm :}
\par
{\rm (1)} For $n= 2m$, 
\begin{equation*}
\CH (Y_n) \otimes \Z/2 \simeq \Z/2 [c_1, c_2] /(b_{m-1}, c_2 b_{m-2})
\otimes \Delta (v_{2m-4}, v_{2m-2}).
\end{equation*}
\par
{\rm (2)} For $n= 2m+1$, 
\begin{equation*}
\CH (Y_n) \otimes \Z/2 \simeq \Z/2 [c_1, c_2] /(b_{m-1}, c_2 b_{m-2})
\otimes \Delta (v_{2m-2}, v_{2m}),
\end{equation*}
where $|c_1|=2, |c_2|=4, |b_i| = 2i$ and $|v_i| = i$. 
\end{thm}

\begin{thm}[\cite{KKT}] Let $p$ be an odd prime. Then we have a ring isomorphism{\rm :}
\par
{\rm (1)} For $n= 2m$, 
\begin{equation*}
\CH (Y_n) \otimes \Z/p \simeq \Z/p [c_1, c_2, \chi_{2m-4}] /
(c_2 \chi_{2m-4}, \chi_{2m-4}^2- d_{m-2}, d_{m-1}), 
\end{equation*}
where $\chi_{2m-4} \in H^{2m-4} (B SO_{2m-4}, \Z/p)$ is the Euler class. 
\par
{\rm (2)} For $n= 2m+1$, 
\begin{equation*}
\CH (Y_n) \otimes \Z/p \simeq \Z/p [c_1, c_2] /(d_{m-1}, c_2^2 d_{m-2}). 
\end{equation*}
\end{thm}

We recall the definitions of $b_i$, $d_i$ and $v_i$. In a polynomial ring $\Z[\alpha, \beta]$,
we set $c_1= \alpha+ \beta$ and $c_2= \alpha \beta$. Then $b_k$ and $d_k$ are defined by
\begin{equation*}
 b_k = (-1)^k \sum_{i=0}^k \alpha^i \beta^{k-i}\\
\end{equation*}
\par
and\\
\begin{equation*}
 d_k = (-1)^k \sum_{i=0}^k \alpha^{2i}\beta^{2k-2i}.
\end{equation*}
\par
The element $v_{2r} \in CH^{2r} (Y_n)$ is defined by 
\par
(1) For $n=2m$, 
\begin{equation}
\begin{cases}
2 v_{2m-4}= \chi_{2m-4}-b_{m-2}\\
2v_{2m-2}= b_{m-1}.
\end{cases}
\end{equation}
\par
(2) For $n=2m+1$, 
\begin{equation}
\begin{cases}
2 v_{2m-2}= b_{m-1}\\
2v_{2m}= c_2 b_{m-2}.
\end{cases}
\end{equation}
\par
The following formulas are easily proved.
\begin{lem} We have 
\begin{align*}
& b_k = (-1)^k \sum_{\mu=0}^{\left[ \frac{k}{2} \right]}
(-1)^\mu {k-\mu \choose \mu} c_1^{k- 2 \mu} c_2^\mu\\
\intertext{and}
&d_k = (-1)^k \sum_{\mu=0}^k
(-1)^\mu {2k-\mu+1 \choose \mu} c_1^{2k- 2 \mu} c_2^\mu.
\end{align*}
\end{lem}
The following lemmas are also easily shown. 
\begin{lem}
\begin{equation*}
\sum_{\mu=0}^h (-1)^\mu c_2^{h-\mu} b_{2 \mu} = d_h.
\end{equation*}
\end{lem}

\begin{lem} We set $f_n (x)= (1+x)^n- (1+x^n)$ 
and write $f_n (x)$ as 
\begin{equation*}
f_n (x)= \sum_{\mu=1}^{\left[ \frac{n}{2} \right]}
a_\mu x^\mu (1+x)^{n-2 \mu}.
\end{equation*}
Then we have 
\begin{equation*}
a_\mu = (-1)^{\mu+1} \frac{n}{\mu} {n-1-\mu \choose \mu-1}.
\end{equation*}
Especially, the last term is given by 
\begin{equation*}
\begin{cases}
(-1)^{s+1} 2 x^s & \quad \text{for $n= 2s$}\\
(-1)^{s+1} (2s+1) x^s (1+x) & \quad \text{for $n= 2s+1$.}
\end{cases}
\end{equation*}
\end{lem}
\par

For $n= 2m$ or $2m+1$, we define a subgroup 
$A_n$ of $\CH (Y_n)$ by 
\begin{equation}
A_n = \left( \underset{i=0}{\overset{m-2}{\oplus}} \Z[c_1]/(c_1^{m-1-i}) 
\{ c_2^i \} \right) \otimes B_n,
\end{equation}
where we set 
\begin{equation*}
B_n = 
\begin{cases}
\Delta_{\Z} (v_{2m-4}, v_{2m-2}) & \quad n=2m\\
\Delta_{\Z} (v_{2m-2}, v_{2m}) & \quad n= 2m+1.
\end{cases}
\end{equation*}
The generators $v_{2i}$ is specified in (3.2) and (3.3). 

The following lemma is proved in the same way as in \cite[Lemma 3.8]{KKT}. 
\begin{lem} For a prime $p$, we abbreviate $\CH (Y_n) \otimes \Z_{(p)}$ 
as $\CH (Y_n)_{(p)}$. 
If $p$ is odd, we have the following isomorphism of modules{\rm :}
\par
{\rm (i)} For $n= 2m$, 
\begin{equation*}
\CH (Y_n)_{(p)} \simeq \Z_{(p)}[c_1]/(c_1^{2(m-1)}) \{1, \chi_{2m-4} \}
\oplus \underset{i=1}{\overset{m-2}{\oplus}} \Z_{(p)}[c_1]/(c_1^{2(m-1-i)}) 
\{ c_2^{2i-1}, c_2^{2i} \}.
\end{equation*}
\par
{\rm (ii)} For $n=2m+1$, 
\begin{equation*}
\CH (Y_n)_{(p)} \simeq \underset{i=0}{\overset{m-2}{\oplus}}
\Z_{(p)} [c_1]/ (c_1^{2(m-1-i)}) \{c_2^{2i}, c_2^{2i+1} \}.
\end{equation*}
\end{lem}

\begin{thm} An integral basis of $\CH (Y_n)$ is constructed from 
the monomial basis of $A_n$ in (3.4). 
The results are summed up in Sect.~5, 5.2.
\end{thm}

The proof is based on a rather complicated calculation. Its outline is as follows: 
We construct a set of suitable generators starting from the basis of $A_n$. 
It is easily verified that it is a basis of $\CH (Y_n) \otimes \Z/p$ 
by using the presentation of Lemma 3.6. Then it is a $\Z$-basis of 
$\CH (Y_n)$. We only prove the case $n= 2m+1$ and $m$ is even. To simplify the proof, we set: 
\begin{align}
& C = \CH (Y_n)\\
& \b_{m-1}= (-1)^{m-1} b_{m-1} = \sum_{i+j=m-1} \alpha^i \beta^j, 
\quad e_m= c_2 \b_{m-2} \notag \\
& \d_{m-1}= (-1)^{m-1} d_{m-1} \notag \\
& \A = \left(\Z[c_1]/(c_1^{m-1}) \oplus \dots \oplus
\Z[c_1]/(c_1^{m-1-i}) \{ c_2^i \} \oplus \dots \oplus \{ c_2^{m-2} \} \right)
\otimes \Z\{1, \b_{m-1}, e_m, \b_{m-1} e_m \}. \notag
\end{align}
Let $\A_n$ the set of a basis 
\begin{equation*}
\{ c_1^i c_2^j \b_{m-1}^{\epsilon_1} e_m^{\epsilon_2}:
i+j \leq m-2, \epsilon_k= 0 \; \text{or}\; 1 \, (k=1, 2)\}.
\end{equation*}
Then we have 
\begin{equation*}
C \subset {\mathbb Q}[c_1, c_2]/(\d_{m-1}, c_2^2 \d_{m-2}), \quad
A \subset {\mathbb Q}[c_1, c_2]/(\d_{m-1}, c_2^2 \d_{m-2})
\end{equation*}
\begin{equation*}
\text{and} \quad 
c_2^{2i} \d_{m-1-i} \in (\d_{m-1}, c_2^2 \d_{m-2}), \; i \geq 0.
\end{equation*}
\par
Our concern is a homogeneous polynomial algebra $S_{\Z} (V)$ of $V=\{c_1, c_2\}$. 
It is identified with an inhomogeneous ring $\Z[x]$ by putting $x= \frac{\beta}{\alpha}$. 
Then we have 
\begin{equation*}
c_1= 1+x, \; c_2=x, \; \b_{m-1} = \frac{x^m-1}{x-1}, \quad \text{and} \quad 
\d_{m-1}= \frac{x^{2m}-1}{x^2-1}.
\end{equation*}
Using the identification, the next formulas are directly checked: In $\Z[c_1, c_2]$, 
\begin{align*}
&-c_2^{2i+1+j} \b_{m-2-2i} = c_2^{j+1} \b_{2i-1} \b_{m-1}-c_2^j \b_{2i} e_m\\
\intertext{and}
&-c_2^{2i+2+j} \b_{m-3-2i} = c_2^{j+1} \b_{2i} \b_{m-1}-c_2^j \b_{2i+1} e_m. 
\end{align*}
We set 
\begin{align}
&\left[ c_1^{2i-1}c_2^{j+1} \b_{m-1} \right]:= c_2^{j+1} \b_{2i-1} \b_{m-1}-c_2^j \b_{2i} e_m = -c_2^{2i+j+1} \b_{m-2-2i}\\
\intertext{and}
&\left[ c_1^{2i}c_2^{j+1} \b_{m-1} \right]:=c_2^{j+1} \b_{2i} \b_{m-1}-c_2^j \b_{2i+1} e_m 
=-c_2^{2i+j+2} \b_{m-3-2i}. \notag
\end{align}
\par
Noting $c_2^2 \d_{m-2}=0$, we have 
\begin{equation*}
c_2^2 \b_{m-2} \b_{m-1}- c_2^2 c_1 \d_{m-2}= \frac{x^3}{(x-1)^2} 
(x^{m-2}-1) (x^{m-1}-1) = c_2^3 \b_{m-3} \b_{m-2}.
\end{equation*}
Using $c_2^{2i} \d_{m-1-i}=0, \; 1 \leq i \leq m$, we repeat the argument and get 
\begin{equation}
c_2^i \b_{m-1} e_m= c_2^{2i+1} \b_{m-2-i} \b_{m-1-i}.
\end{equation}
\par
In $Y_n$, we see that $\d_{m-1}=0$. Recall that we set $f_{2i+3}(x) = 
(1+x)^{2i+3}- (1+x^{2i+3})$ (see Lemma 3.5). Then we have
\begin{align*}
c_1^{2i+1} &c_2^j \b_{m-1} e_m- c_2^{j+1}\d_i \d_{m-1}\\
&=\frac{x^{j+1}(x^m-1)}{(x^2-1)^2} \left((1+x)^{2i+3} (x^{m-1}-1)-(x^{2i+2}-1)(x^m+1)\right)\\
&= \frac{x^{j+1}(x^m-1)}{(x^2-1)^2} \left((1+x^{2i+3}) (x^{m-1}-1)-(x^{2i+2}-1)(x^m+1)+
f_{2i+3}(x) (x^{m-1}-1)\right)\\
&= \frac{x^{2i+j+3}}{(x-1)(x^2-1)}(x^m-1)(x^{m-2i-3}-1)+
\frac{x^{j+1}}{(x^2-1)^2}(x^{m-1}-1)(x^m-1) f_{2i+3} (x).
\end{align*}
We set 
\begin{align*}
I_1 =& \frac{x^{2i+j+3}}{(x-1)(x^2-1)} (x^m-1)(x^{m-2i-3}-1)\\
\intertext{and}
I_2=& \frac{x^{j+1}}{(x^2-1)^2}(x^{m-1}-1)(x^m-1) f_{2i+3}(x).
\end{align*}
Since $x^{2 \lambda} \d_{m-1-\lambda}=0$ for $0 \leq \lambda \leq m-1$, 
we see that $x^{2i+j+3+\lambda} \d_{m-1-(i+2)-\lambda}=0$ for $0 \leq \lambda \leq j-1$. 
Hence
\begin{equation*}
I_1- \sum_{\lambda=0}^{j-1} x^{2i+j+3+\lambda} \d_{m-1-(i+2)-\lambda}=
\frac{x^{2i+2j+3}}{(x-1)(x^2-1)}(x^{m-j}-1)(x^{m-2i-j-3}-1).
\end{equation*}
From Lemma 3.5, 
\begin{equation*}
I_2 = \sum_{\mu=1}^i a_\mu c_2^{j+\mu} c_1^{2i+1-2\mu} \b_{m-1} e_m+
(-1)^i k x^{i+j+2} \frac{(x^{m-1}-1)(x^m-1)}{(x-1) (x^2-1)},
\end{equation*}
where $k= 2i+3$.
\par
We set $J= x^{i+j+2} \frac{(x^{m-1}-1)(x^m-1)}{(x-1)(x^2-1)}$. Using 
$x^{i+j+2+\lambda} \d_{m-2-\lambda}=0$ for $1 \leq \lambda \leq i+j$, we see that 
\begin{equation*}
J- \sum_{\lambda=0}^{i+j} x^{i+j+2+\lambda} \d_{m-2-\lambda}=
\frac{x^{2i+2j+3}}{(x-1)(x^2-1)} (x^{m-2-i-j}-1)(x^{m-1-i-j}-1).
\end{equation*}
Afterwards, we introduce the notation 
\begin{equation}
\left[ c_1^{2i+1}c_2^j \b_{m-1}e_m \right]:= c_1^{2i+1}c_2^j \b_{m-1}e_m-
\sum_{\mu=1}^i a_\mu c_2^{j+\mu}c_1^{2i+1-2\mu} \b_{m-1}e_m,
\end{equation}
where $f_{2i+3}(x) = \sum_{\mu=1}^{i+1} a_\mu x^\mu (1+x)^{2i+3-2\mu}$ 
(see Lemma 3.5). 
\par
We sum up these arguments:
\begin{equation}
\left[c_1^{2i+1} c_2^j \b_{m-1}e_m \right]= I_1+(-1)^i k J,
\end{equation}
where
\begin{align*}
&I_1 = \frac{x^{2i+2j+3}}{(x-1)(x^2-1)} (x^{m-j}-1)(x^{m-2i-j-3}-1),\\
&J=  \frac{x^{2i+2j+3}}{(x-1)(x^2-1)}(x^{m-2-i-j}-1)(x^{m-1-i-j}-1), 
\quad \text{and} \quad k=2i+3.
\end{align*}
\par
A direct calculation shows $J-I_1= \frac{x^{m+j}}{(x-1)(x^2-1)}(x^{i+1}-1)(x^{i+2}-1)$.  
Hence we get a key formula 
\begin{equation}
\left[
c_1^{2i+1} c_2^j \b_{m-1} e_m \right]=
\left( (-1)^i (2i+3)+1 \right) I_1+ (-1)^i (2i+3) \frac{x^{m+j}}{(x-1)(x^2-1)}(x^{i+1}-1)(x^{i+2}-1).
\end{equation}
We note that (3.9) holds for $m$ even or odd. From now on, we assume that $m$ is even. 
When we put $j=1$ in (3.9), we see that $I_1= c_2^{2i+4} \d_{\frac{m}{2}-i-3} e_m$. 
We set 
\begin{equation*}
\langle c_1^{2i+1} c_2 \b_{m-1} e_m \rangle:=
\left[ c_1^{2i+1} c_2 \b_{m-1} e_m \right]-\left((-1)^i(2i+3)+1 \right) c_2^{2i+4} 
\d_{\frac{m}{2}-i-3} e_m.
\end{equation*}
Then we have
\begin{equation}
\frac {\langle c_1^{2i+1} c_2 \b_{m-1} e_m \rangle}{2i+3}=(-1)^i c_2^{m+1} 
\frac{\b_i \b_{i+1}}{c_1}.
\end{equation}
\par
We call a generator $c_1^{m-3-2i}c_2^{2i+1} \b_{m-1} e_m$ to be the head of the 
presentation $\A_n$ in (3.4). For the head generator, it is directly shown that $I_1=0$ by (3.9). Hence 
the formula (3.9) implies
\begin{equation}
\frac{\left[
c_1^{m-3-2i} c_2^{2i+1} \b_{m-1} e_m \right]}{m-1-2i}=
(-1)^i \frac{\b_{\frac{m}{2}-i-2} \b_{\frac{m}{2}-i-1}}{c_1}.
\end{equation}
Then the formulas (3.10) and (3.11) imply the following relations: We set 
\begin{align*}
\langle c_1^{2i+1} \b_{m-1} e_m \rangle:=&
\left[ c_1^{2i+1} \b_{m-1} e_m \right]-\left((-1)^i(2i+3)+1\right) c_2^{2i+2} 
\d_{\frac{m}{2}-i-2} e_m\\
\intertext{and}
\langle c_1^{2 \alpha+1-2 \beta} c_2^{1+ \beta}\b_{m-1} e_m \rangle:=&
\left[ c_1^{2 \alpha+1-2 \beta} c_2^{1+ \beta}\b_{m-1} e_m \right]\\
& \qquad -\left((-1)^{\alpha-\beta}(2\alpha-2\beta+3)+1\right) c_2^{2\alpha+4} 
\d_{\frac{m}{2}-\alpha-3} e_m.
\end{align*}
Then we have 
\begin{align}
&\qquad \langle c_1^{2i+1} \b_{m-1} e_m \rangle - \frac{(-1)^i (2i+3)+1}{(-1)^{i-1}(2i+1)}
\langle c_1^{2i-1}c_2 \b_{m-1} e_m \rangle = -\frac{c_2^m}{c_1}\b_i \b_{i+1}\\
&\qquad \langle c_1^{2\alpha+1-2\beta} c_2^{1+\beta} \b_{m-1} e_m \rangle - 
\frac{(-1)^\beta (2\alpha-2\beta+3)+(-1)^\alpha}{2\alpha+3}
\langle c_1^{2\alpha+1}c_2 \b_{m-1} e_m \rangle\\
&\qquad \qquad = -\frac{c_2^{m+\beta+2}}{c_1}\b_{\alpha-\beta-1} \b_{\alpha-\beta} \notag\\
\intertext{and}
&\qquad \left[ c_1^{m-3-2\alpha-2\beta} c_2^{2\alpha+\beta+1} \right]- 
\frac{(-1)^\beta (m-1-2\alpha+2\beta)+(-1)^{\frac{m}{2}-\alpha}}{m-1-2\alpha}
\left[ c_1^{m-3-2\alpha}c_2^{2\alpha+1} \b_{m-1} e_m \right] \\
&\qquad \qquad= -\frac{c_2^{m+2\alpha+\beta+1}}{c_1}\b_{\gamma} \b_{\gamma+1}, \quad
\text{where $\gamma= \frac{m}{2}-2-\alpha-\beta$.} \notag
\end{align}
\par
Last we consider a generator $c_1^{2i+2}c_2^j \b_{m-1}e_m$. We set 
\begin{align*}
\langle c_1^{2i+2}c_2^j \b_{m-1} e_m \rangle := &c_1^{2i+2}c_2^j \b_{m-1}e_m-
\sum_{\mu=1}^i a_\mu c_2^{j+\mu}c_1^{2i+2-2\mu} \b_{m-1}e_m\\
& \qquad-\left((-1)^i (2i+3)+1\right) c_2^{i+j+1} \b_{m-1}e_m,
\end{align*}
where $f_{2i+3} (x)= \sum_{\mu=1}^{i+1} a_\mu x^\mu (1+x)^{2i+3-2\mu}$ (see Lemma 3.5). 
Then the formula (3.9) implies that 
\begin{equation}
\langle c_1^{2i+1}c_2^j \b_{m-1}e_m \rangle = -c_2^{m+j} \b_i \b_{i+1}. 
\end{equation}
\par
We consider a set 
\begin{align*}
&\{ \text{(3.6), (3.10), (3.11), (3.12), (3.13), (3.14)} \}\\
& \qquad \cup \{x \in \A_n: x = c_1^i c_2^j e_m^{\epsilon_1}, x = c_1^i \b_{m-1},
i+j \leq m-2, \epsilon_k= 0 \; \text{or}\; 1 \}.
\end{align*}
When we reduce it to the mod $p$ reduction for an odd prime $p$, 
the identities from (3.6) to (3.14) show that they are linearly independent 
from Lemma 3.6. Hence, replacing $\b_{m-1}$ and $e_m$ by 
$v_{2m-2}$ and $v_{2m}$ respectively, 
we obtain an integral basis of $\CH (Y_{2m+1})$ for even $m$. 
An integral basis of $\CH (Y_n)$ for other cases is written down in a table of Sect.~5. 
\par

For a group $G$ and a subgroup $H$, let $[G:H]$ denote the index, i.e. the cardinality 
of $G/H$. As a corollary of the above theorem, we have 
\begin{cor} {\rm (i)} For $n= 2m$, 
\begin{equation*}
[\CH (Y_n): A_n] =
\begin{cases}
1^2 \cdot 3^2 \cdot \ldots \cdot (m-3)^2 \cdot (m-1) & \quad \text{$m${\rm :} even}\\
1^2 \cdot 3^2 \cdot \ldots \cdot (m-2)^2 & \quad \text{$m${\rm :} odd.}
\end{cases}
\end{equation*}
\par
{\rm (ii)} For $n= 2m+1$, 
\begin{equation*}
[\CH (Y_n): A_n] =
\begin{cases}
1^2 \cdot 3^2 \cdot \ldots \cdot (m-3)^2 \cdot (m-1) & \quad \text{$m${\rm :} even}\\
1^2 \cdot 3^2 \cdot \ldots \cdot (m-2)^2 \cdot m & \quad \text{$m${\rm :} odd.}
\end{cases}
\end{equation*}
\end{cor}
The following theorem is proved by using the integral basis of $\CH (Y_n)$ given in Theorem 3.7.
\begin{thm} The ring structure of $\CH (Y_n)$ is determined. The results are summed up 
in tables in 5.3 and 5.4 in Sect.~5. (See also 5.5.) 
\end{thm}

\begin{pf} We only show the formula (iii) in Sect. ~5, 5.4. We use the notations of the proof of 
Theorem 3.7. Using $x^{2i} \d_{m-1-i}=0$, we have 
\begin{equation*}
c_1^{m-2i-1} c_2^{2i} \b_{m-1}- x^{2i} \d_{m-1-i}=
\frac{x^m}{x^2-1}(x^{2i}-1)+ \frac{x^{2i}}{x^2-1}(x^m-1) f_{m-2i}(x).
\end{equation*}
We set $I_1= \frac{x^m}{x^2-1} (x^{2i}-1)$ and 
$I_2= \frac{x^{2i}}{x^2-1}(x^m-1) f_{m-2i}(x)$. Then we have 
\begin{equation*}
I_2= \sum_{\mu=1}^{\frac{m}{2}-i-1} a_\mu c_1^{m-2i-1-2\mu} c_2^{2i+\mu}
\b_{m-1}+ (-1)^{\frac{m}{2}-i+1} 2 x^{\frac{m}{2}+i} \frac{x^m-1}{x^2-1}. 
\end{equation*}
We set $J= x^{\frac{m}{2}+i} \frac{x^m-1}{x^2-1}$. From (3.6), we have 
\begin{equation*}
J- c_2^{\frac{m}{2}+i-1} e_m- \left( \sum_{\mu=0}^{\frac{m}{2}-2-i}
(-1)^\mu \left[ c_1^{1+2\mu} c_2^{\frac{m}{2}+i-1-\mu} \b_{m-1} \right] \right)=
(-1)^{\frac{m}{2}-i} x^m \frac{x^{2i}-1}{x^2-1}.
\end{equation*}
Hence
\begin{align}
c_1^{m-2i-1}c_2^{2i} \b_{m-1}-& \sum_{\mu=1}^{\frac{m}{2}-i-1} a_\mu
c_1^{m-2i-1-2\mu} c_2^{2i+\mu} \b_{m-1}\\
&+ (-1)^{\frac{m}{2}-i} 2 \left( c_2^{\frac{m}{2}+i-1} e_m +\sum_{\mu=0}^{\frac{m}{2}-2-i}
(-1)^\mu \left[c_1^{2\mu+1} c_2^{\frac{m}{2}+i-1-\mu} \b_{m-1}\right]\right)=
-x^m \frac{x^{2i}-1}{x^2-1}. \notag
\end{align}
On the other hand, we see from (3.13) that 
\begin{equation}
-x^m \frac{x^{2i}-1}{x^2-1} = \langle c_1^{2i-1} \b_{m-1} e_m \rangle +
\frac{2i+1}{2i-1} \langle c_1^{2i-3} c_2 \b_{m-1} e_m \rangle.
\end{equation}
\par
We calculate 
\begin{equation*}
S:= (-1)^{\frac{m}{2}-i} 2 \left( c_2^{\frac{m}{2}+i-1} e_m +\sum_{\mu=0}^{\frac{m}{2}-2-i}
(-1)^\mu \left[c_1^{2\mu+1} c_2^{\frac{m}{2}+i-1-\mu} \b_{m-1}\right]\right).
\end{equation*}
By (3.6), 
\begin{equation*}
S= (-1)^{\frac{m}{2}-i} 2 (S_1+S_2), 
\end{equation*}
where 
\begin{equation*}
S_1= \sum_{\mu=0}^{\frac{m}{2}-2-i} (-1)^\mu 
c_2^{\frac{m}{2}+i+1-\mu} \b_{2 \mu+1} \b_{m-1} \quad \text{and} \quad 
S_2= c_2^{\frac{m}{2}+i-1} e_m -\sum_{\mu=0}^{\frac{m}{2}-2-i} (-1)^\mu 
c_2^{\frac{m}{2}+i-2-\mu} \b_{2 \mu+2} e_m.
\end{equation*}
We have from Lemma 3.3 that 
\begin{equation*}
S_1= \sum_{\mu=1}^{\frac{m}{2}-1-i} (-1)^{\frac{m}{2}-1-\mu-i} 
{m-1-2i-\mu \choose \mu-1} c_1^{m-2i-1-2\mu} c_2^{2i+\mu} \b_{m-1}.
\end{equation*}
Using $\sum_{\mu=0}^h (-1)^\mu c_2^{h-\mu} \b_{2\mu}= (-1)^h \d_h$ (see Lemma 3.4), 
\begin{equation*}
S_2= c_2^{2i} \left( \sum_{\mu=0}^{\frac{m}{2}-1-i} (-1)^\mu c_2^{\frac{m}{2}-i-1-\mu}
\b_{2 \mu} \right) \b_{m-1}= (-1)^{\frac{m}{2}-1-i} c_2^{2i} \d_{\frac{m}{2}-1-i} 
\b_{m-1}.
\end{equation*}
The formula (3.17) is 
\begin{align*}
& c_1^{m-2i-1}c_2^{2i}\b_{m-1}+ 
\sum_{\mu=1}^{\frac{m}{2}-i-1} (-1)^\mu {m-2i-1-\mu \choose \mu}
c_1^{m-2i-1} c_2^{2i+\mu} \b_{m-1}
-2 c_2^{2i}\d_{\frac{m}{2}-1-i} e_m\\
&=-x^m \frac{x^{2i-1}}{x^2-1}. \notag
\end{align*}
Comparing with (3.18), we obtain
\begin{align*}
&c_1^{m-2i-1}c_2^{2i} \b_{m-1} = \sum_{\mu=1}^{\frac{m}{2}-i-1} (-1)^{\mu+1}
{m-2i-1-\mu \choose \mu} c_1^{m-2i-1} c_2^{2i+\mu} \b_{m-1}\\
& \qquad -\frac{2}{2i-1} c_2^{2i} \d_{\frac{m}{2}-1-i} e_m
+ \left[c_1^{2i-1} \b_{m-1} e_m \right]+ \frac{2i+1}{2i-1}
\left[c_1^{2i-3}c_2 \b_{m-1}e_m \right]. \notag
\end{align*}
This completes the proof of the formula. \qed
\end{pf}

\section{The Chow ring of $X_n$}

\begin{thm} Let $T (X_n)$ and $F(X_n)$ be the torsion part and the free part of 
$\CH (X_n)$, respectively. Then we have 
\par
{\rm (i)} For $n=4t$, 
\begin{align*}
&F(X_n) \simeq \Z[c_2]/(c_2^t) \{1, v_{4t-4}\}\\
&T (X_n) \simeq \Z/2[c_2]/(c_2^t) \{ v_{4t-2}, v_{4t-4} v_{4t-2}\}.
\end{align*}
\par
{\rm (ii)} For $n=4t+1$, 
\begin{align*}
&F(X_n) \simeq \Z[c_2]/(c_2^t)
\oplus \Z[c_2]/(c_2^{t-1})\{v_{4t}\}\\
&T (X_n) \simeq \Z/2[c_2]/(c_2^t) \{ v_{4t-2}, v_{4t-2} v_{4t}\}
\oplus \Z/(2 t) \{c_2^{t-1} v_{4t}\}.
\end{align*}
\par
{\rm (iii)} For $n=4t+2$, 
\begin{align*}
&F(X_n) \simeq \Z[c_2]/(c_2^t) \{1, v_{4t}\}
\oplus \Z\{v_{4t-2}\}\\
&T (X_n) \simeq \Z/2[c_2]/(c_2^{t-1}) \{ c_2 v_{4t-2}, c_2 v_{4t-2} v_{4t}\}
\oplus \Z/4 \{v_{4t-2} v_{4t}\}.
\end{align*}
\par
{\rm (iv)} For $n=4t+3$, 
\begin{align*}
&F(X_n) \simeq \Z[c_2]/(c_2^t) \{1, v_{4t}\}\\ 
&T (X_n) \simeq \Z/2[c_2]/(c_2^t) \{ v_{4t+2}, v_{4t} v_{4t+2}\}
\oplus \Z/(2t+1) \{c_2^t v_{4t} \}.
\end{align*}
\end{thm} 

\begin{pf} Let ${\tilde X}_n=
X_n \times_{\G_m} {\mathbb A}^1$ be the associated bundle of (3.1) and 
$s: Y_n \to {\tilde X}_n$ 
the $0$-section. Since $s^\ast: \CH ({\tilde X}_n) \overset{\sim}{\to} \CH (Y_n)$ by Theorem 2.4 (ii), 
the first assertion of the same theorem for $V= {\tilde X}_n$ and $X= s(Y_n)$ gives an exact sequence 
\begin{equation*}
\CH (Y_n) \overset{\cdot c_1}{\to} \CH (Y_n) \overset{\pi^\ast}{\to} 
\CH (X_n) \to 0. 
\end{equation*}
Theorem 4.1 follows from this and the ring structure of $\CH (Y_n)$ 
in Theorems 3.7 and 3.9. \qed
\end{pf}

Next we consider the cycle map. The cohomology groups mean an etale 
cohomology \cite{G}, \cite{M}. 
All varieties are defined over $K'$, which is a subfield of an 
algebraically closed field $K$. Let $l$ be a prime with $(l, {\rm ch} (K))=1$. 
We denote a locally constant sheaf $\mu_l^{\otimes i}$ by $\Z/l^{(i)}$. 

\begin{cor} The homomorphism $cl: CH^i (X_n) \to H^{2i} (X_n, \Z/l^{(i)})$ is 
injective. 
\end{cor}
\begin{pf} Since $({\tilde X}_n, Y_n)$ is a smooth pair, we have the Gysin sequence 
as in \cite[Appendice 1.3.3]{Del2} and \cite[VI Remark 5.4]{M}. Since the cycle map and the Gysin map are commutative, 
we have the following commutative diagram, where each row is exact: 
\begin{equation*}
\begin{CD}
CH^i (Y_n) @> \cdot c_1 >> CH^{i+1} (Y_n) @> \pi^\ast>> CH^{i+1} (X_n) @>>> 0\\
@VV cl V @VV cl V @VV cl V\\
H^{2i} (Y_n, \Z/l^{(i)}) @> \cdot c_1>> H^{2(i+1)} (Y_n, \Z/l^{(i+1)}) @> \pi^\ast>>
H^{2(i+1)} (X_n, \Z/l^{(i+1)})
\end{CD}
\end{equation*}
Corollary 4.2 follows from Theorem 4.1 and this diagram. \qed
\end{pf}

\begin{rem} Assume that we have a $K'$-isomorphism 
$Y_n \simeq SO (n, K)/(SO (n-4,K) \times 
GL (2, K))$, where $SO (n, K)$ and $SO (n-4, K)$ are split over $K'$, 
and that $\pi: Y_n \to X_n$ is a $K'$-map. Then the Galois actions ${\overline G} =
{\rm Gal} ({\overline K}'/K')$ on $H^\cdot (Y_n, \Z/l^{(i)})$ and 
$H^\cdot (X_n, \Z/l^{(i)})$ are described by the character group $X(T)$ of a 
$K'$-split maximal torus $T$ of $SO (n,K)$. It follows from a result of 
\cite[8-2]{Del2}.
\end{rem}

\section{Tables of the ring structure of $\CH (Y_n)$}
\par
\noindent
{\bf 5.1. Notations.} (i) For $k \in {\mathbb N} \cup \{0\}$, 
we define $b_k$ and $d_k \in {\mathbb Z}[c_1, c_2]$ as follows: 
\begin{align*}
& b_k = (-1)^k \sum_{\mu=0}^{\left[ \frac{k}{2} \right]}
(-1)^\mu {k-\mu \choose \mu} c_1^{k- 2 \mu} c_2^\mu\\
\intertext{and}
&d_k = (-1)^k \sum_{\mu=0}^k
(-1)^\mu {2k-\mu+1 \choose \mu} c_1^{2k- 2 \mu} c_2^\mu.
\end{align*}
\par
(ii) For $g \in {\mathbb N}$ and 
$\mu \in {\mathbb N} \cup \{0, -1\}$, we define $a_{g,\mu} \in {\mathbb Z}$ by 
\begin{equation*}
a_{g,\mu}=
\begin{cases} 
(-1)^{1+\mu} \frac{g}{\mu} {g-1-\mu \choose \mu-1} & \quad \mu \geq 1\\
-1 & \quad \mu=0\\
0 & \quad \mu =-1. 
\end{cases}
\end{equation*}
\par
Then the integers $a_{g, \mu}$ are characterized by
\begin{equation*}
(1+x)^g = 1+ x^g+ \sum_{\mu=1}^{\left[ \frac{g}{2} \right]} a_{g,\mu} x^\mu (1+x)^{g-2\mu}.
\end{equation*}
\par
(iii) The generators $v_{2i}$ are given by (3.2) and (3.3).
\par
\medskip
\noindent
{\bf 5.2. An integral basis of $\CH (Y_n)$.} In the following (I) and (II), 
we give an integral basis of $\CH (Y_n)$. The notations are explained as follows: 
Let $S_n$ be the set of the monomial basis of $A_n$ in (3.4). Let $T$ be a subset of $S_n$. 
Then for an element $\xi \in T$, $\langle \xi \rangle$ (resp. $\langle \xi \rangle'$) 
is defined to be the right-hand side of an equation (1)-(8) below. We consider a set
\begin{equation*}
\left\{ \frac{\langle \xi \rangle}{l_\xi}: \xi \in T \right\} 
\cup 
\left\{ \eta: \eta \in S_n-T \right\}, 
\end{equation*}
where $l_\xi \in {\mathbb N}$. 
Following this procedure, we obtain an integral basis of $\CH (Y_n)$. 
We abbreviate this basis as 
$\displaystyle{
\left\{ \frac{\langle \xi \rangle}{l_\xi}: \xi \in T \right\}.}$
\par
\medskip
(I) The case $n=2m$. 
\par
\quad 
(i) For even $m$, 
\begin{equation*}
\left\{
\frac{\langle c_1^{2i+1} c_2 v_{2m-4} v_{2m-2} \rangle}{2i+3}, 
\frac{\langle c_1^{m-2j-3} c_2^{2j+1} v_{2m-4} v_{2m-2} \rangle}{m-2j-1}:
0 \leq i \leq \frac{m}{2}-2, 1 \leq j \leq \frac{m}{2}-2 \right\}.
\end{equation*}
\par
\quad
(ii) For odd $m$, 
\begin{equation*}
\left\{
\frac{\langle c_1^{2i+1} c_2 v_{2m-4} v_{2m-2} \rangle'}{2i+3}, 
\frac{\langle c_1^{m-2j-2} c_2^{2j} v_{2m-4} v_{2m-2} \rangle'}{m-2j}:
0 \leq i \leq \frac{m-5}{2}, 1 \leq j \leq \frac{m-3}{2} \right\}.
\end{equation*}
\par
\medskip
(II) The case $n= 2m+1$. 
\par
\quad 
(iii) For even $m$, 
\begin{equation*}
\left\{
\frac{\langle c_1^{2i+1} c_2 v_{2m-2} v_{2m} \rangle}{2i+3}, 
\frac{\langle c_1^{m-2j-3} c_2^{2j+1} v_{2m-2} v_{2m} \rangle}{m-2j-1}:
0 \leq i \leq \frac{m}{2}-2, 1 \leq j \leq \frac{m}{2}-2 \right\}.
\end{equation*}
\par
\quad
(iv) For odd $m$, 
\begin{equation*}
\left\{
\frac{\langle c_1^{2i+1} v_{2m-2} v_{2m} \rangle}{2i+3}, 
\frac{\langle c_1^{m-2j-2} c_2^{2j} v_{2m-2} v_{2m} \rangle}{m-2j}:
0 \leq i \leq \frac{m-3}{2}, 1 \leq j \leq \frac{m-3}{2} \right\}.
\end{equation*}
\par
\bigskip
\noindent
Here $\langle \quad \rangle$ and $\langle \quad \rangle'$ are defined as follows: 
{\allowdisplaybreaks
\begin{align*}
&\langle c_1^{2i+1} c_2 v_{2m-4} v_{2m-2} \rangle = c_1^{2i+1} c_2 v_{2m-4} v_{2m-2} 
+ (-1)^{\frac{m+2i+2}{2}} \frac{(-1)^i (2i+3)+1}{2} c_2^{2i+4} d_{\frac{m-2i-6}{2}} 
v_{2m-4} \tag{1}\\
&\phantom{\langle c_1^{2i+1} c_2 v_{2m-4} v_{2m-2} \rangle =}
- \sum_{\mu=1}^i a_{2i+3, \mu} c_1^{2i+1-2\mu} c_2^{1+\mu} v_{2m-4} v_{2m-2}.\\
&\langle c_1^{m-2j-3} c_2^{2j+1} v_{2m-4} v_{2m-2} \rangle = 
c_1^{m-2j-3} c_2^{2j+1} v_{2m-4} v_{2m-2} \tag{2}\\
& \phantom{\langle c_1^{m-2j-3} c_2^{2j+1} v_{2m-4} v_{2m-2} \rangle =}
- \sum_{\mu=1}^{\frac{m-2j-4}{2}} a_{m-2j-1,\mu} c_1^{m-2j-3-2\mu} c_2^{2j+1+\mu} 
v_{2m-4} v_{2m-2}.\\
&\langle c_1^{2i+1} c_2 v_{2m-4} v_{2m-2} \rangle' =  c_1^{2i+1} c_2 v_{2m-4}v_{2m-2}
+(-1)^{\frac{m+2i+1}{2}} \frac{(-1)^i (2i+3)+1}{2} c_2^{2i+3} d_{\frac{m-2i-5}{2}}
v_{2m-2} \tag{3}\\
& \phantom{\langle c_1^{2i+1} c_2 v_{2m-4} v_{2m-2} \rangle' =}
-\sum_{\mu=1}^i a_{2i+3,\mu} c_1^{2i+1-2\mu} c_2^{1+\mu} v_{2m-4} v_{2m-2}.\\
&\langle c_1^{m-2j-2} c_2^{2j} v_{2m-4}v_{2m-2} \rangle'= 
c_1^{m-2j-2} c_2^{2j} v_{2m-4} v_{2m-2} \tag{4}\\
& \phantom{\langle c_1^{m-2j-2} c_2^{2j} v_{2m-4}v_{2m-2} \rangle'=}
-\sum_{\mu=1}^{\frac{m-2j-3}{2}} a_{m-2j,\mu} c_1^{m-2j-2-2\mu}c_2^{2j+\mu}
v_{2m-4} v_{2m-2}.\\
&\langle c_1^{2i+1} c_2 v_{2m-2} v_{2m} \rangle = 
c_1^{2i+1} c_2 v_{2m-2} v_{2m} + (-1)^{\frac{m+2i+2}{2}} \frac{(-1)^i (2i+3)+1}{2}
c_2^{2i+4} d_{\frac{m-2i-6}{2}} v_{2m} \tag{5}\\
& \phantom{\langle c_1^{2i+1} c_2 v_{2m-2} v_{2m} \rangle =}
-\sum_{\mu=1}^i a_{2i+3,\mu} c_1^{2i+1-2\mu}c_2^{1+\mu} v_{2m-2} v_{2m}.\\
&\langle c_1^{m-2j-3} c_2^{2j+1} v_{2m-2} v_{2m} \rangle =
c_1^{m-2j-3} c_2^{2j+1} v_{2m-2} v_{2m} \tag{6}\\
& \phantom{\langle c_1^{m-2j-3} c_2^{2j+1} v_{2m-2} v_{2m} \rangle =}
- \sum_{\mu=1}^{\frac{m-2j-4}{2}} a_{m-2j-1,\mu} c_1^{m-2j-3-2\mu} c_2^{2j+1+\mu} 
v_{2m-2} v_{2m}.\\
&\langle c_1^{2i+1} v_{2m-2} v_{2m} \rangle = c_1^{2i+1} v_{2m-2}v_{2m}+
(-1)^{\frac{m+2i+3}{2}} \frac{(-1)^i (2i+3)+1}{2} c_2^{2i+3} d_{\frac{m-2i-5}{2}}
v_{2m-2} \tag{7}\\
& \phantom{\langle c_1^{2i+1} v_{2m-2} v_{2m} \rangle =}
- \sum_{\mu=1}^i a_{2i+3,\mu} c_1^{2i+1-2\mu} c_2^\mu v_{2m-2} v_{2m}.\\
&\langle c_1^{m-2j-2} c_2^{2j}v_{2m-2} v_{2m} \rangle = 
c_1^{m-2j-2} c_2^{2j} v_{2m-2} v_{2m}
\tag{8}\\
& \phantom{\langle c_1^{m-2j-2} c_2^{2j}v_{2m-2} v_{2m} \rangle =}
- \sum_{\mu=1}^{\frac{m-2j-3}{2}} a_{m-2j, \mu} c_1^{m-2j-2-2\mu} c_2^{2j+\mu} 
v_{2m-2} v_{2m}.
\end{align*}}
\par
\medskip
\noindent
{\bf 5.3. The ring structure of $\CH (Y_n)_{(2)}$ for $n=2m$.} 
\begin{table}[H]
\begin{center}
\begin{tabular}{c|c|c}
\noalign{\hrule height0.8pt}
& even $m$ & odd $m$\\
\noalign{\hrule height0.8pt}
$c_1^{m-1}$  &
\multicolumn{2}{c} {$(1)$}\\
\hline 
$c_1^{m-k-1} c_2^k \, (k \geq 1)$ & 
\multicolumn{2}{c}{$(2)$}\\
\noalign{\hrule height0.8pt}
$c_1^{m-1} v_{2m-4}$ &
\multicolumn{2}{c}{$(3)$}\\
\hline
$c_1^{m-2i-1} c_2^{2i} v_{2m-4} \, (i \geq 1)$ & 
\multicolumn{2}{c}{$(4)$}\\
\hline
$c_1^{m-2i-2} c_2^{2i+1} v_{2m-4} \, (i \geq 0)$ & $(5)$ & $(6)$\\
\noalign{\hrule height0.8pt}
$c_1^{m-2i-1} c_2^{2i} v_{2m-2} \, (i \geq 0)$ & 
$(7)$ & $(8)$\\
\hline
$c_1^{m-2i-2} c_2^{2i+1} v_{2m-2} \, (i \geq 0)$ & 
\multicolumn{2}{c}{$(9)$}\\
\noalign{\hrule height0.8pt}
$c_1^{m-2i-1} c_2^{2i} v_{2m-4} v_{2m-2} \, (i \geq 0)$ & 
$(10)$ & $(11)$\\
\hline
$c_1^{m-2i-2} c_2^{2i+1} v_{2m-4} v_{2m-2} \, (i \geq 0)$ 
& $(12)$ & $(13)$\\
\noalign{\hrule height0.8pt}
$v_{2m-4}^2$ & 
$(14)$ & $(15)$\\
\hline
$v_{2m-2}^2$ 
& $(16)$ & $(17)$\\
\noalign{\hrule height0.8pt}
\end{tabular}
\end{center}
\end{table}
\par
\noindent
Here
{\allowdisplaybreaks
\begin{align*}
(1) = & 
\left\{ \sum_{\mu=1}^{\left[ \frac{m-1}{2} \right]} 
(-1)^{1+\mu} \binom{m-1-\mu}{\mu} 
c_1^{m-1-2\mu} c_2^\mu \right\}+ (-1)^{m+1}2 v_{2m-2}.\\
(2) = & \left\{ \sum_{\mu=1}^{\left[ \frac{m-k-1}{2} \right]} 
(-1)^{1+\mu} \binom{m-k-1-\mu}{\mu} 
c_1^{m-k-1-2\mu} c_2^{k+\mu} \right\} \\
&+ \left\{ (-1)^{m+k} 2c_2 b_{k-1} \right\} v_{2m-4} + 
\left\{ (-1)^{m+k}2 c_2 b_{k-2} \right\} v_{2m-2}.\\
(3) = & \left\{ \sum_{\mu=1}^{\left[ \frac{m-1}{2} \right]}
(-1)^{1+\mu} \binom{m-1-\mu}{\mu} c_1^{m-1-2\mu} c_2^\mu \right\} v_{2m-4}+ 
(-1)^{m+1} 2 v_{2m-4} v_{2m-2}.\\
(4) = & \left\{ \sum_{\mu=1}^{\left[ \frac{m-2i-1}{2} \right]} 
(-1)^{1+\mu} \binom{m-2i-1-\mu}{\mu}
c_1^{m-2i-1-2\mu} c_2^{2i+\mu} \right\} v_{2m-4}\\
& + \left\{ (-1)^m 2 \sum_{\mu=0}^{i-1} a_{2i-1, \mu} c_1^{2i-2-2\mu} c_2^{1+\mu} \right\}
v_{2m-4} v_{2m-2}.\\
(5) = & \left\{ (-1)^{\frac{m+2i+2}{2}} \frac{4i}{2i+1} c_2^{2i+2} 
d_{\frac{m-2i-4}{2}}
+ \sum_{\mu=1}^{\frac{m-2i-2}{2}} a_{m-2i-1, \mu} c_1^{m-2i-2-2\mu} 
c_2^{2i+1+\mu} \right\} v_{2m-4} \\
&+ \left\{ 
2\sum_{\mu=-1}^{i-2} \left( a_{2i-1, \mu} + \frac{2i-1}{2i+1} a_{2i+1, 1+\mu} \right) c_1^{2i-3-2\mu} c_2^{2+\mu} \right\} v_{2m-4} v_{2m-2}.\\
(6) = & \left\{ \sum_{\mu=1}^{\frac{m-2i-3}{2}} 
(-1)^{1+\mu} \binom{m-2i-2-\mu}{\mu} 
c_1^{m-2i-2-2\mu} c_2^{2i+1+\mu} \right\} v_{2m-4}\\
&+ \left\{ (-1)^{\frac{m+2i+1}{2}} \frac{2}{2i+1} c_2^{2i+1} 
d_{\frac{m-2i-3}{2}} \right\} v_{2m-2}\\
&- \left\{ 
2\sum_{\mu=-1}^{i-2} \left( a_{2i-1, \mu} + \frac{2i-1}{2i+1} a_{2i+1, 1+\mu} \right) c_1^{2i-3-2\mu} c_2^{2+\mu} \right\} v_{2m-4} v_{2m-2}.\\
(7) = & \left\{ (-1)^{\frac{m+2i}{2}} \frac{2}{2i+1} c_2^{2i+2} 
d_{\frac{m-2i-4}{2}} \right\} v_{2m-4}\\
&+ \left\{ \sum_{\mu=1}^{\frac{m-2i-2}{2}} 
(-1)^{1+\mu} \binom{m-2i-1-\mu}{\mu} c_1^{m-2i-1-2\mu} 
c_2^{2i+\mu} \right\} v_{2m-2}\\
&+ \left\{ 
2\sum_{\mu=-1}^{i-2} \left( a_{2i-1, \mu} + \frac{2i-1}{2i+1} a_{2i+1, 1+\mu} \right) c_1^{2i-3-2\mu} c_2^{2+\mu} \right\} v_{2m-4} v_{2m-2}.\\
(8) = & \left\{  (-1)^{\frac{m+2i+3}{2}} \frac{4i}{2i+1} c_2^{2i+1} 
d_{\frac{m-2i-3}{2}} 
+ \sum_{\mu=1}^{\frac{m-2i-1}{2}} 
a_{m-2i, \mu} c_1^{m-2i-1-2\mu} 
c_2^{2i+\mu} \right\} v_{2m-2}\\
&- \left\{ 
2\sum_{\mu=-1}^{i-2} \left( a_{2i-1, \mu} + \frac{2i-1}{2i+1} a_{2i+1, 1+\mu} \right) c_1^{2i-3-2\mu} c_2^{2+\mu} \right\} v_{2m-4} v_{2m-2}.\\
(9) = & \left\{ \sum_{\mu=1}^{\left[ \frac{m-2i-2}{2} \right]} 
(-1)^{1+\mu} \binom{m-2i-2-\mu}{\mu} 
c_1^{m-2i-2-2\mu} c_2^{2i+1+\mu} \right\} v_{2m-2}\\
&+ \left\{(-1)^{m} 
2\sum_{\mu=0}^{i} a_{2i+1, \mu} c_1^{2i-2\mu} 
c_2^{1+\mu} \right\} v_{2m-4} v_{2m-2}.\\
(10) = & \left\{  \sum_{\mu=0}^{\frac{m-2i-4}{2}} 
\left( \frac{m-2i+1}{m-2i-1} a_{m-2i-1, \mu}+ a_{m-2i+1, 1+\mu} \right) 
c_1^{m-2i-3-2\mu} c_2^{2i+1+\mu} \right\} v_{2m-4} v_{2m-2}.\\
(11) = & \left\{ \sum_{\mu=1}^{\frac{m-2i-1}{2}} 
a_{m-2i, \mu}
c_1^{m-2i-1-2\mu} c_2^{2i+\mu} \right\} v_{2m-4} v_{2m-2}.\\
(12) = & \left\{ \sum_{\mu=1}^{\frac{m-2i-2}{2}} 
a_{m-2i-1, \mu}
c_1^{m-2i-2-2\mu} c_2^{2i+1+\mu} \right\} v_{2m-4} v_{2m-2}.\\
(13) = & \left\{  \sum_{\mu=0}^{\frac{m-2i-5}{2}} 
\left( \frac{m-2i}{m-2i-2} a_{m-2i-2, \mu}+ a_{m-2i, 1+\mu} \right) 
c_1^{m-2i-4-2\mu} c_2^{2i+2+\mu} \right\} v_{2m-4} v_{2m-2}.\\
(14) = & \; (-1)^{\frac{m}{2}} d_{\frac{m-2}{2}} v_{2m-4}.\\
(15) = & \; -b_{m-2}  v_{2m-4}+ (-1)^{\frac{m+3}{2}} d_{\frac{m-3}{2}} v_{2m-2}.\\
(16) = & \; (-1)^{\frac{m+2}{2}} c_2^2 d_{\frac{m-4}{2}} v_{2m-4}.\\
(17) = & \; (-1)^{\frac{m+1}{2}} c_2 d_{\frac{m-3}{2}} v_{2m-2}.
\end{align*}}
\par
\newpage
\noindent
{\bf 5.4. The ring structure of $\CH (Y_n)_{(2)}$ for $n=2m+1$.}
\begin{table}[H]
\begin{center}
\begin{tabular}{c|c|c}
\noalign{\hrule height0.8pt}
& even $m$ & odd $m$\\
\noalign{\hrule height0.8pt}
$c_1^{m-1}$  &
\multicolumn{2}{c} {$\text{\rm (i)}$}\\
\hline 
$c_1^{m-k-1} c_2^k \, (k \geq 1)$ & 
\multicolumn{2}{c}{$\text{\rm (ii)}$}\\
\noalign{\hrule height0.8pt}
$c_1^{m-2i-1} c_2^{2i} v_{2m-2} \, (i \geq 0)$ & 
$\text{\rm (iii)}$ & $\text{\rm (iv)}$\\
\hline
$c_1^{m-2i-2} c_2^{2i+1} v_{2m-2} \, (i \geq 0)$ & 
\multicolumn{2}{c} {$\text{\rm (v)}$}\\
\noalign{\hrule height0.8pt}
$c_1^{m-1} v_{2m}$ &
\multicolumn{2}{c}{$\text{\rm (vi)}$}\\
\hline
$c_1^{m-2i-1} c_2^{2i} v_{2m} \, (i \geq 1)$ & 
\multicolumn{2}{c}{$\text{\rm (vii)}$}\\
\hline
$c_1^{m-2i-2} c_2^{2i+1} v_{2m} \, (i \geq 0)$ & $\text{\rm (viii)}$ & 
$\text{\rm (ix)}$\\
\noalign{\hrule height0.8pt}
$c_1^{m-2i-1} c_2^{2i} v_{2m-2} v_{2m} \, (i \geq 0)$ & 
$\text{\rm (x)}$ & $\text{\rm (xi)}$\\
\hline
$c_1^{m-2i-2} c_2^{2i+1} v_{2m-2} v_{2m} \, (i \geq 0)$ 
& $\text{\rm (xii)}$ & $\text{\rm (xiii)}$\\
\noalign{\hrule height0.8pt}
$v_{2m-2}^2$ & 
$\text{\rm (xiv)}$ & $\text{\rm (xv)}$\\
\hline
$v_{2m}^2$ 
& $\text{\rm (xvi)}$ & $\text{\rm (xvii)}$\\
\noalign{\hrule height0.8pt}
\end{tabular}
\end{center}
\end{table}
\par
\noindent
Here
{\allowdisplaybreaks
\begin{align*}
\text{\rm (i)} = & 
\left\{ \sum_{\mu=1}^{\left[ \frac{m-1}{2} \right]} 
(-1)^{1+\mu} \binom{m-1-\mu}{\mu} 
c_1^{m-1-2\mu} c_2^\mu \right\}+ (-1)^{m+1} 2 v_{2m-2}.\\
\text{\rm (ii)} = & \left\{ \sum_{\mu=1}^{\left[ \frac{m-k-1}{2} \right]} 
(-1)^{1+\mu} \binom{m-k-1-\mu}{\mu} 
c_1^{m-k-1-2\mu} c_2^{k+\mu} \right\} \\
&+ \left\{ (-1)^{m+k} 2c_2 b_{k-2} \right\} v_{2m-2} + 
\left\{ (-1)^{m+k+1} 2b_{k-1} \right\} v_{2m}.\\
\text{\rm (iii)} = & \left\{ \sum_{\mu=1}^{\frac{m-2i-2}{2}}
(-1)^{1+\mu} \binom{m-2i-1-\mu}{\mu}  
c_1^{m-2i-1-2\mu} c_2^{2i+\mu} \right\} 
v_{2m-2}\\
& + \left\{ (-1)^{\frac{m+2i+2}{2}} \frac{2}{2i-1} c_2^{2i} 
d_{\frac{m-2i-2}{2}} \right\} v_{2m}\\
&- \left\{ 
2\sum_{\mu=0}^{i-1} \left( \frac{2i+1}{2i-1} a_{2i-1, -1+\mu}+
a_{2i+1, \mu} \right) c_1^{2i-1-2\mu} c_2^\mu \right\} v_{2m-2} v_{2m}.\\
\text{\rm (iv)} = & \left\{ (-1)^{\frac{m+2i+3}{2}} \frac{4i}{2i+1} 
c_2^{2i+1} d_{\frac{m-2i-3}{2}} 
+ \sum_{\mu=1}^{\frac{m-2i-1}{2}} 
a_{m-2i, \mu} c_1^{m-2i-1-2\mu} 
c_2^{2i+\mu} \right\} v_{2m-2}\\
&+ \left\{ 
2\sum_{\mu=-1}^{i-2} \left( a_{2i-1, \mu} + \frac{2i-1}{2i+1} a_{2i+1, 1+\mu} \right) c_1^{2i-3-2\mu} c_2^{1+\mu} \right\} v_{2m-2} v_{2m}.\\
\text{\rm (v)} = & \left\{ \sum_{\mu=1}^{\left[ \frac{m-2i-2}{2} \right]} 
(-1)^{1+\mu} \binom{m-2i-2-\mu}{\mu}
c_1^{m-2i-2-2\mu} c_2^{2i+1+\mu} \right\} v_{2m-2}\\
& + \left\{ (-1)^{m+1}2 \sum_{\mu=0}^{i} a_{2i+1, \mu} c_1^{2i-2\mu} c_2^{\mu} \right\}
v_{2m-2} v_{2m}.\\
\text{\rm (vi)} = & \left\{ \sum_{\mu=1}^{\left[\frac{m-1}{2} \right]} 
(-1)^{1+\mu} \binom{m-1-\mu}{\mu} 
c_1^{m-1-2\mu} c_2^{\mu} \right\} v_{2m} + (-1)^{m+1} 2v_{2m-2} v_{2m}.\\
\text{\rm (vii)} = & \left\{ \sum_{\mu=1}^{\left[ \frac{m-2i-1}{2} \right]} 
(-1)^{1+\mu} \binom{m-2i-1-\mu}{\mu} 
c_1^{m-2i-1-2\mu} c_2^{2i+\mu} \right\} v_{2m}\\
&+ \left\{(-1)^m 2\sum_{\mu=0}^{i-1} a_{2i-1,\mu} c_1^{2i-2-2\mu} c_2^{1+\mu}
\right\} v_{2m-2} v_{2m}.\\
\text{\rm (viii)} =  & \left\{ 
 (-1)^{\frac{m+2i+2}{2}} \frac{4i}{2i+1} c_2^{2i+2} 
d_{\frac{m-2i-4}{2}} +
\sum_{\mu=1}^{\frac{m-2i-2}{2}} a_{m-2i-1, \mu} 
c_1^{m-2i-2-2\mu} c_2^{2i+1+\mu} \right\} v_{2m}\\
&+ \left\{ 
2\sum_{\mu=-1}^{i-2} \left( a_{2i-1, \mu} + \frac{2i-1}{2i+1} a_{2i+1, 1+\mu} \right) c_1^{2i-3-2\mu} c_2^{2+\mu} \right\} v_{2m-2} v_{2m}.\\
\text{\rm (ix)} = & \left\{ (-1)^{\frac{m+2i+3}{2}} 
\frac{2}{2i+3} c_2^{2i+3} 
d_{\frac{m-2i-5}{2}} \right\} v_{2m-2}\\
&+ \left\{ \sum_{\mu=1}^{\frac{m-2i-3}{2}} 
(-1)^{1+\mu} \binom{m-2i-2-\mu}{\mu} c_1^{m-2i-2-2\mu} 
c_2^{2i+1+\mu} \right\} v_{2m}\\
&+ \left\{2 \sum_{\mu=-1}^{i-1} \left(
a_{2i+1, \mu}- a_{2i+1, 1+\mu} + \frac{2i+1}{2i+3} a_{2i+3, 1+\mu} 
\right) c_1^{2i-1-2\mu} c_2^{1+\mu} \right\} v_{2m-2} v_{2m}.\\
\text{\rm (x)} = & \left\{  \sum_{\mu=0}^{\frac{m-2i-4}{2}} 
\left( \frac{m-2i+1}{m-2i-1} a_{m-2i-1, \mu}+ a_{m-2i+1, 1+\mu} \right) 
c_1^{m-2i-3-2\mu} c_2^{2i+1+\mu} \right\} v_{2m-2} v_{2m}.\\
\text{\rm (xi)} = & \left\{ \sum_{\mu=1}^{\frac{m-2i-1}{2}} 
a_{m-2i, \mu}
c_1^{m-2i-1-2\mu} c_2^{2i+\mu} \right\} v_{2m-2} v_{2m}.\\
\text{\rm (xii)} = & \left\{ \sum_{\mu=1}^{\frac{m-2i-2}{2}} 
a_{m-2i-1, \mu} c_1^{m-2i-2-2\mu} c_2^{2i+1+\mu} \right\} v_{2m-2} v_{2m}.\\
\text{\rm (xiii)} = & \left\{  \sum_{\mu=0}^{\frac{m-2i-5}{2}} 
\left( \frac{m-2i}{m-2i-2} a_{m-2i-2, \mu}+ a_{m-2i, 1+\mu} \right) 
c_1^{m-2i-4-2\mu} c_2^{2i+2+\mu} \right\} v_{2m-2} v_{2m}.\\
\text{\rm (xiv)} = & \; (-1)^{\frac{m+2}{2}} d_{\frac{m-2}{2}} v_{2m}.\\
\text{\rm (xv)} = & \; (-1)^{\frac{m+1}{2}} c_2 d_{\frac{m-3}{2}} v_{2m-2}.\\
\text{\rm (xvi)} = & \; (-1)^{\frac{m}{2}} c_2^2 d_{\frac{m-4}{2}} v_{2m}.\\
\text{\rm (xvii)} = & \; (-1)^{\frac{m+1}{2}} \frac{1}{3} c_2^3 d_{\frac{m-5}{2}} v_{2m-2}
- \frac{2}{3} c_1 v_{2m-2} v_{2m}.
\end{align*}}
\par
\medskip
\noindent
{\bf 5.5. A remark on the ring structure of $\CH (Y_n)$.} We have given the ring 
structure of $\CH (Y_n)_{(2)}$ in 5.3 and 5.4. But actually, it is easy to determine 
the ring structure of $\CH (Y_n)$ from 5.2, 5.3 and 5.4. For example, 
by the basis in 5.2, the formula 5.3 (5) is rewritten as follows: 
{\allowdisplaybreaks
\begin{align*}
(5)' \quad c_1^{m-2i-2} c_2^{2i+1} & v_{2m-4} = \Biggl\{ (-1)^{\frac{m+2i+2}{2}} 
\left( (-1)^i (2i-1) +1 \right) c_2^{2i+2} d_{\frac{m-2i-4}{2}}\\
& + \sum_{\mu=1}^{\frac{m-2i-2}{2}} a_{m-2i-1, \mu} c_1^{m-2i-2-2\mu} 
c_2^{2i+1+\mu} \Biggr\} v_{2m-4}\\
&+ \Biggl\{ 2\sum_{\mu=-1}^{i-2} a_{2i-1, \mu} c_1^{2i-3-2 \mu} c_2^{2+\mu} 
\Biggr\} v_{2m-4} v_{2m-2}
- \frac{4i-2}{2i+1} \langle c_1^{2i-1} c_2 v_{2m-4} v_{2m-2} \rangle. 
\end{align*}}
The other cases can be calculated similarly.

\end{document}